\documentclass{article}
\usepackage[utf8]{inputenc}
\usepackage{amsmath,amscd,amssymb,amsfonts,amsthm,courier,relsize,bm}
\usepackage{hyperref,enumerate,mathrsfs,mathtools,slashed}
\usepackage{tikz-cd}
\tikzcdset{ampersand replacement=\&}
\usepackage{xcolor}  
\hypersetup{
    colorlinks,
    linkcolor={red!50!black},
    citecolor={blue!70!black},
    urlcolor={blue!80!black}
}

\newcommand\norm[1]{\left\lVert#1\right\rVert}

\newtheorem{theorem}{Theorem}[section]

\newtheorem{proposition}[theorem]{Proposition}

\theoremstyle{definition}    
\newtheorem{definition}[theorem]{Definition}
\theoremstyle{remark}

\newcommand{\ignore}[1]{}

\def\Spin{\ensuremath{\textnormal{Spin}}}
\def\SO{\ensuremath{\textnormal{SO}}}

\def\bC{\ensuremath{\mathbb{C}}}
\def\bR{\ensuremath{\mathbb{R}}}
\def\bZ{\ensuremath{\mathbb{Z}}}

\def\End{\ensuremath{\textnormal{End}}}

\def\ker{\ensuremath{\textnormal{ker}}}

\def\Tr{\ensuremath{\textnormal{Tr}}}

\def\index{\ensuremath{\textnormal{index}}}

\newcommand{\br}[1]{\bm{\mathrm{#1}}}
\title{The Geometry of Mehler's Kernel}
\author{Jesus Sanchez Jr}
\date{December 2023}

\begin{document}
\maketitle
\begin{abstract}
    We study the relationship between the Getzler calculus of a spin Riemannian manifold and the Riemannian geometry of the corresponding principal spin bundle. We then use the calculus of Gaussian-Grassmann integrals developed by Berline--Vergne to compute the Getzler symbol of the spinor heat flow. 
\end{abstract}
\section{Introduction}
Given a closed smooth manifold $M$ and an elliptic pseudo-differential operator $P:\Gamma(M,E)\rightarrow\Gamma(M,F)$ acting between the sections of complex vector bundles over $M$, global analysis tells us that $P$ is Fredholm and the Atiyah-Singer index theorem \cite{Ati-Sin} computes the Fredholm index via a topological expression
\[
    \text{Index}P=(-1)^{\text{dim} M}\int_{T^*M}ch(\sigma(P))\wedge \widehat{Td}(TM^{\bC}).
\]
To date, there are many approaches to the index theorem in the general case of pseudo-differential operators, most of which focus on the global aspects of $\sigma(P)$ as an element of a (possibly noncommutative) generalized cohomology theory, e.g. $K$-theory. While these approaches are incredibly powerful and have built various bridges with other subjects in mathematics such as operator $K$-theory, deformation quantization, and Lie groupoids to name a few, these approaches also wash away the relationship between the operator $P$ and geometric structures on $M$. The subject of local index theory is aimed at deriving index theorems by studying the interplay of the local coefficients of $P$ and a suitably chosen geometric structure $(M, G)$ on $M$.\par 
The case when $P=D$ is a Dirac-type operator and the geometric structure is Riemannian has been well studied in the literature and begins with the work of Atiyah-Bott-Patodi \cite{atiyah1973heat}. The major drawback of \cite{atiyah1973heat} was its implicit extraction of the characteristic forms relevant for the index calculation. This was later clarified with the work of Alvarez-Guamé \cite{AG} and Witten \cite{Wit} on supersymmetric quantum mechanics and localization on loop space. Rigorous proofs of these claims have been given, the three techniques which stand out are of Bismut \cite{Bis1}, Berline--Vergne \cite{BV1}, and Getzler \cite{Get}, each technique using and developing different machinery to handle the problem. The goal of the present article is to build a bridge between the two different approaches of Getzler \cite{Get} and Berline--Vergne \cite{BV1}. We will show that one can use the Berline--Vergne approach to evaluate the Getzler symbol of the heat kernel, demystifying the geometry of the terms within the Getzler-Mehler kernel, in particular shedding light on the geometry of the Gaussian term. In all approaches to the local index theorem, the heat semi-group $\{e^{-t\slashed{D}^2}\}_{t>0}$ plays a vital role. How one approaches the small time behaviour of $e^{-t\slashed{D}^2}$ is different for Getzler and Berline--Vergne. Let us briefly explain the approaches.\par
In Getzler's original approach, one views the super-manifold $C^{\infty}(\mathscr{N})=\Gamma(T^*M,\pi^*\wedge T^*M)$ as the cotangent superbundle bundle to the super-manifold $\mathcal{M}$ with $C^{\infty}(\mathcal{M})=\Gamma(M,\mathscr{S})$, the sections of the spinor bundle. While this idea is heuristic (because $\Gamma(M,\mathscr{S})$ is not a true supermanifold) it motivates the definition of the Getzler symbol calculus as a combination of the dequantization of the pseudodifferential operators on $C^{\infty}(M)$ to functions on $T^*M$ and the dequantization of the Clifford algebra $\Gamma(M,C\ell(T^*M))$ to the exterior algebra $\Gamma(M,\wedge T^*M)$. In a follow up article \cite{Get2}, Getzler explains the symbol calculus as a rescaling at the level of Schwartz kernel which is particularly simple for differential operators and smoothing operators. The main achievement of this technique is to show that the symbol of the heat kernel is the Mehler kernel 
\begin{equation*}
    \sigma^{Par-G}(e^{-t\slashed{D}^2})(\br{x})=(4\pi t)^{-n/2}\text{det}^{1/2}\left( \frac{tR/2}{\text{sinh}(tR/2)}\right)\text{exp}\left( -\frac{1}{4t}\left\langle \br{x}\middle|\frac{tR}{2}\text{coth}\left( \frac{tR}{2}\right)\middle|\br{x}\right\rangle\right) 
\end{equation*}
which, when evaluated at $(\br{x},t)=(\br{0},1)$, recovers the quantity 
\begin{equation*}
(4\pi)^{-n/2}\text{det}^{1/2}\left(\frac{R/2}{\text{sinh}(R/2)}\right),
\end{equation*}
the $\hat{A}$-genus of the manifold $M$.\par 
Coming at the problem from a different angle, Berline--Vergne studied the geometry of a principal spin bundle $P_{Spin}(M)$ of a given Riemannian manifold $(M,g)$. The key reason for doing this is that one can lift $\slashed{D}^2$ on the base $M$ to $P_{Spin}(M)$ via
\begin{displaymath} 
\begin{tikzcd}[column sep=8em]
C^{\infty}\left(P_{Spin}(M),\wedge\mathbb{C}^{n/2}\right)^{\rho} \arrow[r, "\Delta+\pi^*k/4-n(n-1)/2"] \arrow[d, "\simeq"]
\& C^{\infty}\left(P_{Spin}(M),\wedge\mathbb{C}^{n/2}\right)^{\rho} \arrow[d, "\simeq" ] \\
\Gamma(\mathscr{S}) \arrow[r, "\slashed{D}^2"]
\& \Gamma(\mathscr{S})
\end{tikzcd}
\end{displaymath}
which gives rise to the following expression for the heat kernel on the base
\begin{equation}
\label{EA}
\left\langle p\middle|e^{-t\slashed{D}^2}\middle|q\right\rangle=e^{n(n-1)t/2}\int_{\Spin(n)}\left\langle p\cdot g\middle|e^{-t(\Delta+\pi^*k/4)}\middle|q\right\rangle\rho(g)dg.
\end{equation}
The appearance of the $\hat{A}$-genus of $M$ can be understood by making two observations. The first is that, from this perspective, it is only the \emph{first} term in the heat kernel asymptotics of $\langle p|e^{-t(\Delta+\pi^*k/4)}|q\rangle$ which contributes to the local supertrace in the limit as $t\rightarrow 0^+$. The second is that the first term in the heat asymptotics of $\langle p|e^{-t(\Delta+\pi^*k/4)}|q\rangle$ only depends on the local Jacobian of the exponential map. It is here that the geometry of $P_{Spin}(M)$ is utilized to show that the $\hat{A}$-genus appears \emph{explicitly} by studying the Jacobian of the exponential map restricted to the horizontal distribution defined by the Levi--Civita connection on $M$. 
\begin{theorem}[Berline--Vergne, \cite{BV1}]
Fix $p\in P_{Spin}(M)$, $\pi(p)=x$ and consider the function $\Phi$ defined for $A\in\mathfrak{spin}(n)$ small
\begin{align*}
    \Phi(A)&= \textup{det}^{-1/2}(\textup{exp}_p)_{*,A}\\
    &=\textup{det}^{-1/2}\left(\frac{\sinh ad_A/2}{ad_A/2}\right)\cdot\textup{det}^{-1/2}\left(\frac{\sinh\tau(A\cdot\Omega_p)/4}{\tau(A\cdot\Omega_p)/4}\right).
    \end{align*}
Then one has (after making an identification $\bigwedge \mathbb{R}^n\simeq\bigwedge T^*_xM$)
\begin{align*}
    \Phi(2e\wedge e^*)&=\textup{det}^{1/2}\left(\frac{R_x/2}{\sinh (R_x/2)}\right).\\
\end{align*}
\end{theorem}
On the surface, these two methods seem unrelated. Moreover, Getzler does more by constructing a full symbol calculus and recovering the full Mehler kernel, not just the $\hat{A}$-genus. In this article, we will use the techniques developed by Berline--Vergne to uncover the geometry behind the Mehler kernel. Here is the approach: following Berline--Vergne we will use the equivariant averages \eqref{EA} and heat kernel asymptotics of the operator $\Delta+\pi^*k/4$ on $P_{Spin}(M)$ to obtain an asymptotic expansion of $e^{-t\slashed{D}^2}$. What we will find is that the terms of the heat kernel asymptotic expansion are expressed as Gaussian averages along the fibers of $P_{Spin}(M)$. From here we will combine the Riemannian geometry of $P_{Spin}(M)$ and the Getzler rescaling on the base $M$ to lift the Getzler rescaling to $P_{Spin}(M)$ (see Sec 4.2). The combination of the Getzler rescaling and the Gaussian averages will convert the problem to an evaluation of Gaussian-Grassman integrals. The latter have already been developed by Berline--Vergne, our only addition will be to work with Taylor expansions of such integrals. The added insight is that the Gaussian term of the Mehler kernel is governed by the geodesic distance function on $P_{Spin}(M)$.
\begin{theorem}
Fix $p\in P_{Spin}(M)$, $\pi(p)=x$, and $v\in\mathcal{H}_p$ small. Consider the function $\Psi$ defined for $A\in\mathfrak{spin}(n)$ small
\begin{align*}
    \Psi(A)&= \left\langle v,\left(\frac{I-\textup{exp}(-\tau(A\cdot\Omega_p)/2)}{\tau(A\cdot\Omega_p)/2}\right)^{-1}v\right\rangle.
    \end{align*}
Then one has (after making an identification $\bigwedge \mathbb{R}^n\simeq\bigwedge T^*_xM$)
\begin{align*}
    \Psi(2e\wedge e^*)&=\left\langle v, \frac{R_x}{2}\coth\left(\frac{R_x}{2}\right)v\right\rangle
\end{align*}
\end{theorem}
When studying the asymptotics of the Getzler-rescaled heat kernel through the lenses of the geometry on $P_{Spin}(M)$, we fully recover the Mehler kernel. 
\begin{theorem}
Fix $p\in P_{Spin}(M)$, $\pi(p)=x$, and $v\in\mathcal{H}_p$ small. For $0<|u|\leq 1$, we have 
\begin{align*}
    &\sigma(k_{u^2}(\exp_p(uv),p))\\
    &\sim (4\pi u)^{-n}\sum_{j=0}^{n/2}u^{2j}\left(\textup{det}^{1/2}\left(\frac{R_x/2}{\sinh(R_x/2)}\right)\textup{exp}\left( -\frac{1}{4}\left\langle v,\frac{R_x}{2}\coth\left( \frac{R_x}{2}\right)v\right\rangle\right)\right)_{2j}\\
    &\hspace{1cm}+ u^{-n+1}\sum_{j=0}^{n/2}u^{2j}\omega_{2j}(\vec{v})+\cdots
\end{align*}
where $\omega_{2j}(v)$ is analytic in $v$. Evaluating $\lim_{u\rightarrow 0}u^n\sum_{j=0}^nu^{-j}\sigma_j(k_{u^2}(\exp_p(uv),p))$ gives
\[
(4\pi )^{-n/2}\textup{det}^{1/2}\left(\frac{R_x/2}{\sinh(R_x/2)}\right)\exp\left( -\frac{1}{4}\left\langle v,\frac{R_x}{2}\coth\left( \frac{R_x}{2}\right)v\right\rangle\right).
\] 
\end{theorem}
\subsection*{Acknowledgments}
The author is indebted to Matthew Bernstein and Sergio Zamora Barrera who gave many insightful comments to the first draft of the article. The author is further indebted to Ahmad Reza Haj Saeedi Sadegh for the exciting discussions on this topic which influenced this work.

\section{The Setup}
We briefly orient the reader with our set of notations and background theorems we will be using throughout the article. We consider $(M^n,g)$ to be a closed, oriented, spin Riemannian manifold of even dimension $n$. We choose an equivariant double cover $P_{Spin}(M)$ of $P_{SO}(M)$ given by the diagram below (R is the right action)
\begin{displaymath} 
\begin{tikzcd}[column sep=8em]
P_{Spin}(M)\times \Spin(n) \arrow[r, "R"] \arrow[d, "p\times \widetilde{Ad}"]
\& P_{Spin}(M) \arrow[d, "p" ] \\
P_{SO}(M)\times \SO(n) \arrow[r, "R"]
\& P_{SO}(M)
\end{tikzcd}
\end{displaymath}
If we let $\omega^{LC}\in\Omega^1_{Ad,v}(P_{SO}(M),\mathfrak{so}(n))$ denote the Levi-Civita connection 1-form, then we obtain the induced spinorial Levi-Civita connection 1-form $\omega^{LC}_s:=\tau^{-1}\left(p^*\omega^{LC}\right)\in\Omega^1_{Ad,v}(P_{Spin}(M),\mathfrak{spin}(n))$ where $\tau=D(\widetilde{Ad}):\mathfrak{spin}(n)\xrightarrow{\sim}\mathfrak{so}(n)$. If $\rho:\Spin(n)\rightarrow U(\wedge\mathbb{C}^{n/2})$ denotes the complex unitary spinor representation, then we form the complex vector bundle of spinors over $M$ via $\mathscr{S}:=P_{Spin}(M)\times_{\rho}\wedge\mathbb{C}^{n/2}$. Since the spinor representation $\rho$ is unitary, the spinor bundle $\mathscr{S}$ comes equipped with a Hermitian fiber metric. The spinor bundle comes equiped with the spinorial Levi-Civita connection, denoted $\nabla^{LC}$, which is defined by the commutative diagram
\begin{displaymath} 
\begin{tikzcd}[column sep=8em]
\Omega_{\rho,h}^k(P_{Spin}(M),\wedge\mathbb{C}^{n/2}) \arrow[r, "d+\rho\left(\omega^{LC}_s\right)\wedge"] \arrow[d, "\simeq"]
\& \Omega^{k+1}_{\rho,h}(P_{Spin}(M),\wedge\mathbb{C}^{n/2}) \arrow[d, "\simeq" ] \\
\Omega^k(M,\mathscr{S}) \arrow[r, "d^{\nabla^{LC}}"]
\& \Omega^{k+1}(M,\mathscr{S})
\end{tikzcd}
\end{displaymath}
The Levi-Civita connection on $\mathscr{S}$ is compatible with the Hermitian fiber metric.
The spinor bundle comes with an action of the Clifford bundle of $T^*M$,
\[c:\Gamma(\mathbb{C}\ell(T^*M,g))\times\Gamma(\mathscr{S})\rightarrow\Gamma(\mathscr{S})
\]
By simplicity and dimension comparison we have that $\mathbb{C}\ell(T^*M,g)\simeq_c\End{\mathscr{S}}$ where the isomorphism is of algebra-bundles over $M$. The Clifford action has the property that $c(\omega)^*=-c(\omega)$, for every $\omega\in\Gamma(T^*M)$, where the $^*$ denotes the adjoint operation in $\Gamma(\End{\mathscr{S}})$ given by the Hermitian fiber metric. This allows us to define the main operator of interest for us, the spinor Dirac operator $\slashed{D}$, given by the following diagram 
\begin{displaymath}
  \begin{tikzcd}[column sep=2em]
   \Gamma (\mathscr{S}) \ar[dr,"\slashed{D}"'] \arrow{r}{\nabla^{LC} } \&  \Omega ^1(M,\mathscr{S})  \arrow{d}{\text{c} }\\
  \& \Gamma (\mathscr{S}).
  \end{tikzcd}
\end{displaymath}
The orientation on $M$ guarantees the existence of the volume form $dV\in\Gamma(\wedge^nT^*M)$. We define the complex volume form $dV_{\mathbb{C}}:=i^{n/2}dV$, which has the property that $c(dV_{\mathbb{C}})^2=1$. This allows us to introduce a $\mathbb{Z}_2$-grading on the spinor bundle via $\mathscr{S}^{\pm}:=(1\pm c(dV_{\mathbb{C}}))\mathscr{S}$. All of the structure introduced allows us to conclude that $\slashed{D}$ is an odd, self-adjoint, first order elliptic operator on $\Gamma(\mathscr{S})$. In particular, the operator 
\[\slashed{D}^+:=\slashed{D}|_{\Gamma(\mathscr{S}^+)}:\Gamma(\mathscr{S}^+)\longrightarrow \Gamma(\mathscr{S}^-)
\]
is a first order elliptic operator between the sections of $\mathscr{S}^+$
and $\mathscr{S}^-$. We wish to compute its Fredholm index. 
\subsection{Heat Flow and Index Theory}
The classic route to computing $\index{ \, \slashed{D}^+}$ is to study the interplay between the local coefficients of the operator $\slashed{D}$ and the geometry of $(M,g)$. The first step is to study heat flow on spinors. Functional analytic properties of $\slashed{D}$ guarantee the existence of the heat semi-group $\{e^{-t\slashed{D}^2}\}_{t\geq0}$ as a semi-group of bounded operators acting on $L^2(M,\mathscr{S})$. Further Sobolev theory guarantees the existence of a smooth time dependent kernel $k_t\in\Gamma(\mathscr{S}\boxtimes\mathscr{S}^*)$ for $t>0$, such that 
\begin{itemize}
    \item $e^{-t\slashed{D}^2}s(x)=\int_Mk_t(x,y)s(y)\text{d}V(y)$,\hspace{.5cm}$s\in\Gamma(\mathscr{S})$
    \item $(\partial_t+\slashed{D}^2_x)k_t(x,y)\equiv 0$ for all $(t,x,y)\in  (0,\infty ) \times M\times M$.
\end{itemize}
Thus $e^{-t\slashed{D}^2}$ is a smoothing operator, and we call $k_t$ the heat kernel. The $\pm$-grading on $\mathscr{S}$ induces a $\pm$-grading on $L^2(M,\mathscr{S})$ and thus a $\pm$-grading on $B(L^2(M,\mathscr{}{S}))$. We use this to define a graded trace $\text{Str}:\mathscr{L}^1(L^2(M,\mathscr{S}))\rightarrow\bC$, via the formula
\[\text{Str}(T):=\Tr(c(dV_{\bC})T),\hspace{.5cm} T\in\mathscr{L}^1(L^2(M,\mathscr{S})).
\]
The reason for studying heat flow in this setting is because of the Mckean-Singer formula \cite{MS}, which computes $\text{Index}{\, \slashed{D}^+}$ as
\[\text{Index}{\, \slashed{D}^+}=\text{Str}(e^{-t\slashed{D}^2}),\hspace{.5cm} t>0.
\]
Using the heat kernel, we conclude that 
\begin{align*}
    \text{Index}{\, \slashed{D}^+}&=\text{Str}(e^{-t\slashed{D}^2})\\
        &=\int_M\text{str}(k_t(x,x))dV(x).
\end{align*}
where $\text{str}:\Gamma(\End{\mathscr{S}})\rightarrow C^{\infty}(M)$ is the local supertrace given by
\[\text{str}(T):=\text{tr}(c(dV_{\bC})T),\hspace{.5cm} T\in\Gamma(\End\mathscr{S}).
\]
The goal of \textbf{\emph{local index theory}} is to understand the small time behaviour of $\text{str}(k_t(x,x))$. The reason for doing this is as follows. Functional analytic arguments show that as $t\rightarrow\infty$, the heat kernel approaches the smooth kernel given by projection onto the kernel of $\slashed{D}$ and in particular completely delocalizes (by delocalize we mean that when evaluating the kernel on a section, the result depends on the values of the section all across the manifold, rather than only local values of the given section). Thus we conclude that in the limit $t\rightarrow\infty$, the kernel $k_t(x,y)$ shouldn't contain any local geometric information. However, in the limit $t\rightarrow 0^+$, functional analytic arguments show that the kernel $k_t(x,y)$ localizes around the diagonal of $M\times M$. 
Moreover, intuition tells us that the small time behaviour of heat flow should be governed by the local geometry of our manifold $M$ and the bundle $\mathscr{S}$. 
\subsection{Differential Forms and The Local Supertrace}
As we mentioned earlier, there needs to be a bridge between the local coefficients of our operator and the local geometry of the manifold. In the spinorial setting it is given by the isomorphism $\Gamma(\bC\ell(T^*M,g))\simeq\Gamma(\End\mathscr{S})$. More precisely, we have that $\bC\ell(T^*M,g)$ is isomorphic, as a vector bundle, to $\wedge T^*M^{\bC}$, which is given by
\[\sigma:\bC\ell(T^*M,g)\xrightarrow{\sim}\wedge T^*M^{\bC},\hspace{.5cm}\sigma(e^{i_1}\cdots e^{i_j})=e^{i_1}\wedge\cdots\wedge e^{i_j},
\]
where $\{e^i\}_{i=1}^n$ is a local orthonormal frame for $T^*M$. This means that we can view a section $T\in\Gamma(\End\mathscr{S})$ as Clifford multiplication by a differential form on the base $M$. Part of the utility in this description of $\Gamma(\End\mathscr{S})$ is in its ability to express $\text{str}$ in geometric terms. To see this, first define the \textbf{\emph{Berezin integral}} $\mathcal{B}:\Gamma(\wedge T^*M^{\bC})\rightarrow C^{\infty}(M)$, which is
\[\mathcal{B}(\omega)=\omega_n,
\]
where $\omega_ndV$ is the component of $\omega$ in $\wedge^nT^*M^{\bC}$. Since $\wedge T^*M^{\bC}$ is a bundle of graded commutative algebras, the Berezin integral is a graded trace. Moreover, given $T\in\Gamma(\End\mathscr{S})\simeq\Gamma(\bC\ell(T^*M,g))$, we have
\[\text{str}(T)=(2/i)^{n/2}\mathcal{B}(\sigma(T)).
\]
\section{A Recap of Getzler's Work}
For our purposes, the topic in Getzler's work we want to focus on is the \textbf{\emph{Getzler rescaling}} appearing in \cite{Get2}. The key will be to use the  differential form perspective on $\Gamma(\End\mathscr{S})$ to define a symbol calculus for $\Gamma(\mathscr{S}\boxtimes\mathscr{S}^*)$ which will allow us to understand the local behaviour of $\kappa\in\Gamma(\mathscr{S}\boxtimes\mathscr{S}^*)$ around the diagonal in $M\times M$.\par
\begin{definition}
\label{GR}
Consider a smooth kernel $\kappa\in\Gamma(\mathscr{S}\boxtimes\mathscr{S}^*)$, a fixed point $y\in M$, and a normal neighborhood $(U_y,\vec{x})$ of $y$. For a given $x\in U_y$, let $P_{x\rightarrow y}:\mathscr{S}_x\xrightarrow{\sim}\mathscr{S}_y$ denote parallel translation along the unit speed radial geodesic connecting $x$ and $y$ in $U_y$. We obtain
\[\boldsymbol{\kappa}_y(\vec{x}):=P_{x\rightarrow y}(\kappa(x,y))\in C^{\infty}(U_y,\End\mathscr{S}_y),\hspace{.5cm} x=\exp^M_y(\vec{x})
\]
Using $\End\mathscr{S}_y\simeq_c\bC\ell(T^*_yM)\simeq_{\sigma}\wedge T^*_yM^{\bC}$, we get
\[\sigma(\boldsymbol{\kappa}_y(\vec{x}))\in C^{\infty}(U_y,\wedge T^*_yM^{\bC})
\]
The \textbf{\emph{Getzler rescaling}} of $\kappa$ at $y\in M$ is defined to be 
\[\delta_{u,y}(\kappa)(\vec{x}):=\sum_{j=0}^nu^{-j}\sigma_j(\boldsymbol{\kappa}_y(u\vec{x})),\hspace{.5cm} 0<|u|\leq 1.
\]
\end{definition}
Note that if we allow $u$ to vary then the Getzler rescaling of $\kappa$ at $y$ may be viewed as a smooth function
\[
    \delta_y(\kappa):[-1,1]^{\times}\times U_y\rightarrow\wedge T^*_yM^{\bC},
\]
where $[-1,1]^{\times}$ denotes $[-1,1]$ without $0$.
\begin{definition}
We say that $\kappa\in\Gamma(\mathscr{S}\boxtimes\mathscr{S}^*)$ has Getzler order $m\in\bZ$ at $y\in M$ if $\delta_y(\kappa)$ extends past $u=0$ to a smooth function $\delta_y(\kappa):[-1,1]\times U_y\rightarrow \wedge T^*_yM^{\bC}$. In this case,
the limit below 
\[\sigma^G_{m,y}(\kappa)(\vec{x}):=\lim_{u\rightarrow 0}u^m\delta_{u,y}(\kappa)(\vec{x})\in S(T^*M)\otimes\wedge T^*_yM^{\bC}.
\]
defines the \textbf{\emph{Getzler symbol}} of $\kappa$ at $y\in M$. Elements of $S(T^*M)\otimes\wedge T^*_y M^{\bC}$ are to be viewed as polynomials on $T_yM$ with coefficients in the algebra $\wedge T^*_yM^{\bC}$.
\end{definition}
The way this computation is done in practice is to Taylor expand $\delta_{u,y}(\kappa)(\vec{x})$ in $\vec{x}$ and collect terms of the same $u$ power. When this is done, an aymptotic expansion of $\delta_{u,y}(\kappa)(\vec{x})$ is obtained that has the form
\begin{align*}
    \delta_{u,y}(\kappa)(\vec{x})&\sim\sum_{j=-n}^{\infty}u^j\gamma_{-j}(\vec{x})\\
    &\sim u^{m}\gamma_{-m}(\vec{x})+u^{m+1}\gamma_{-m-1}(\vec{x})+\cdots, \hspace{.5cm} m\geq -n
\end{align*}
where each $\gamma_j(\vec{x})$ is a polynomial in $\vec{x}$ whose coefficients lie in the algebra $\wedge T^*_yM^{\bC}$. The Getzler symbol of $\kappa$ at $y$ will be the first $m$ such that $\gamma_{-m}\neq 0$,
\[\sigma^G_{-m,y}(\kappa)(\vec{x})=\gamma_{-m}(\vec{x})
\]
It is clear from its definition that the Getzler symbol of $\kappa$ gives information about the local behavior of $\kappa$ around the diagonal of $M\times M$. Moreover, we have a relation for computing local supertraces. If the Getzler order of $\kappa$ at y is $m$, then
\begin{equation*}
\text{str}(\kappa(y,y))=\begin{cases}
          0 \quad &\text{if} \,\,\, m<n \\
          (2/i)^{n/2}\mathcal{B}\big(\sigma^G_{n,y}(\kappa)(\vec{0})\big) \quad &\text{if} \,\,\, m=n \\
     \end{cases}
\end{equation*}
A slight modification of this rescaling will be needed to deal with time dependent sections of $\mathscr{S}\boxtimes\mathscr{S}^*$.
\begin{definition}
Suppose $\kappa_t$ is a smooth time dependent section of $\mathscr{S}\boxtimes\mathscr{S}^*$, $t>0$. Repeating the construction of \eqref{GR}, we obtain a smooth time dependent function
\[\sigma(\boldsymbol{\kappa}_{t,y})(\vec{x})\in C^{\infty}((0,\infty)\times U_y,\wedge T^*_yM^{\bC})
\]
The \textbf{\emph{parabolic Getzler rescaling}} of $\kappa_t$ at $y$ to be
\[\delta_{u,y}(\kappa_t)(\vec{x}):=\sum_{j=0}^nu^{-j}\sigma_j(\boldsymbol{\kappa}_{u^2t,y}(u\vec{x})),\hspace{.5cm} 0<|u|\leq 1.
\]
We say in this case that $\kappa_t$ has \textbf{\emph{parabolic Getzler order}} $m\in\bZ$ at $y\in M$ if $u^m\delta_{u,y}(\kappa_t)$ extends past $u=0$ to a smooth function $u^m\delta_{u,y}(\kappa_t):[-1,1]\times(0,\infty)\times U_y\rightarrow \wedge T^*_yM^{\bC}$ and define the \textbf{\emph{parabolic Getzler symbol}} of $\kappa_t$ via the limit below
\[
    \sigma^{Par-G}_{m,y}(\kappa_t)(\vec{x}):=\lim_{u\rightarrow 0}u^m\delta_{u,y}(\kappa_t)(\vec{x})\in C^{\infty}((0,\infty)\times U_y,\wedge T^*_yM^{\bC}).
\]
\end{definition}
It is not always the case that if we are given $\kappa_t$, its parabolic Getzler rescaling $\delta_{u,y}(\kappa_t)(\vec{x})$ will have good behaviour as $u\rightarrow 0$. To take an example, consider multiplying \emph{any} $\kappa\in\Gamma(\mathscr{S}\boxtimes\mathscr{S}^*)$ by $e^{-1/t}$. The time dependent section $e^{-1/t}\kappa$ will not have good rescaled behaviour in the limit as $u\rightarrow 0$. Part of the reason for this is that as $t\rightarrow 0$, we do not have any amount of control on the behaviour of a general $\kappa_t$, i.e. we cannot Taylor expand at $t=0$. However, it will be the case for us that the rescaled heat kernel has great behaviour as $u \rightarrow 0$. We aim to understand this behaviour, and our main technique will be to follow Berline--Vergne and lift the problem to $P_{Spin}(M)$.
\section{A Recap of Berline-Vergne's Approach}
In the approach of Berline--Vergne \cite{BV1}, one should lift the heat kernel of $e^{-t\slashed{D}^2}$ to $P_{Spin}(M)$. In doing so, one is allowed to use the extra symmetries present in $P_{Spin}(M)$ to enable a simple approach to the small time asymptotics of $e^{-t\slashed{D}^2}$. Let $\mathcal{V}$ denote the vertical subbundle of $TP_{Spin}(M)$. If we consider $\langle\cdot,\cdot\rangle :\mathfrak{spin}(n)\rightarrow\bR$ to be an Ad-$\Spin(n)$-invariant inner product on $\mathfrak{spin}(n)$, then using the global trivialization of $\mathcal{V}\simeq P_{Spin}(M)\times\mathfrak{spin}(n)$ allows us to transfer the inner product $\langle\cdot,\cdot\rangle$ to a fiber metric on $\mathcal{V}$. Furthermore, if we let $\mathcal{H}:=\ker\,\omega^{LC}_s$, then $\pi^*TM\simeq\mathcal{H}$ and we use this isomorphism to pull the metric $g$ on $M$ to a fiber metric on $\mathcal{H}$. The splitting $TP_{Spin}(M)=\mathcal{V}\oplus\mathcal{H}$ defines a Riemannian geometry on $P_{Spin}(M)$ by declaring the sum $\mathcal{V}\oplus\mathcal{H}$ to be orthogonal. We note that this geometry on $P_{Spin}(M)$ is acted on isometrically by the right action of $\Spin(n)$. Now let $\Delta^S$ denote the scalar (or vector valued) Riemannian Laplacian on $P_{Spin}(M)$. 
\begin{proposition}[\cite{BV1}]
\label{ComDia}
We have the following commutative diagram
\begin{displaymath} 
\begin{tikzcd}[column sep=8em]
C^{\infty}(P_{Spin}(M),\wedge\mathbb{C}^{n/2})^{\rho} \arrow[r, "\Delta^S+\pi^*k/4-n(n-1)/2"] \arrow[d, "\simeq"]
\& C^{\infty}(P_{Spin}(M),\wedge\mathbb{C}^{n/2})^{\rho} \arrow[d, "\simeq" ] \\
\Gamma(\mathscr{S}) \arrow[r, "\slashed{D}^2"]
\& \Gamma(\mathscr{S})
\end{tikzcd}
\end{displaymath}
\end{proposition}
Recall that the product representation $\rho\times\rho^t:\Spin(n)\times \Spin(n)\rightarrow \End \wedge\bC^{n/2} $ given by
\[(\rho\times\rho^t)(g,h)T=\rho(g)T\rho(h)^{-1},\hspace{.5cm} (g,h)\in \Spin(n)\times \Spin(n),\hspace{.1cm} T\in\End \wedge\bC^{n/2} 
\]
gives rise to the following isomorphism of vector bundles \[\Gamma(\mathscr{S}\boxtimes\mathscr{S}^*)\simeq C^{\infty}(P_{Spin}(M)\times P_{Spin}(M),\End \wedge\bC^{n/2})^{\rho\times\rho^t}.
\]
Thus, for each $t>0$, we may consider the heat kernel $k_t$ as a smooth function on $P_{Spin}(M)\times P_{Spin}(M)$ with values in $\End\wedge\bC^{n/2}$ which is equivariant with respect to the $\Spin(n)\times \Spin(n)$ action. Using \eqref{ComDia}, we obtain the following equality relating the heat kernel $h_t$ of $e^{-t(\Delta^S+\pi^*k/4)}$ to $k_t$
\begin{equation}
\label{EqH}
k_t(p,q)=e^{tn(n-1)/2}\int_{Spin(n)}h_{t}(p\cdot g,q)\rho(g)dg,\hspace{.5cm} (p,q)\in P_{Spin}(M)\times P_{Spin}(M).
\end{equation}
\subsection{Heat Kernel Asymptotics}
The operator $\Delta^S+\pi^*k/4$ is a generalized Laplacian which acts on scalar (or vector) valued functions on $P_{Spin}(M)$. Following \cite{BGV}, there is a formal power series 
\[ (4\pi t)^{-n(n+1)/4}e^{-d(p,q)^2/4t}\sum_{j=0}^{\infty}t^j\Phi_j(p,q)
\]
with $\Phi_j$ smooth scalar functions on $P_{Spin}(M)\times P_{Spin}(M)$ satisfying the following properties:
\begin{itemize}
    \item[(I)] $d(p,q)$ is the Riemannian distance on $P_{Spin}(M)$
    \item[(II)] There is a small neighborhood $U_{\Delta}$ of the diagonal such that each $\Phi_j$ is supported in this neighborhood
    \item[(III)] In a neighborhood $V_{\Delta}\subset U_{\Delta}$ of the diagonal, the formal power series formally solves the heat equation
    \[(\partial_t+\Delta^S_p+\pi^*k(p)/4)(4\pi t)^{-n(n+1)/4}e^{-d(p,q)^2/4t}\sum_{j=0}^{\infty}t^j\Phi_j(p,q)=0
    \]
    \item[(IV)] In the neighborhood $V_{\Delta}$, the $\Phi_j$ are given by
    \begin{align*}
        \Phi_0(\exp^S_p(\vec{x}),p)&=\text{det}^{-1/2}((\exp^S_p)_{*,\vec{x}})\\
        \Phi_j(\exp^S_p(\vec{x}),p)&=-\text{det}^{-1/2}((\exp^S_p)_{*,\vec{x}})\int_0^1t^{j-1}\left(\Delta^S_p+\pi^*k(p)/4\right)\Phi_{j-1}(\exp^S_p(\vec{tx}),p)dt
    \end{align*}
\end{itemize}
Lastly, the formal power series above approximates the heat kernel $h_t$ in the following sense:
\begin{theorem}[\cite{BGV} Ch.2, \cite{Roe} Ch.7]
\label{hka}
For every $j,\ell\geq0$ there exists $k_j$ such that for every $k\geq k_j$ there is a $C_{j,k,\ell}>0$ such that
\[\norm{h_t-(4\pi t)^{-n(n+1)/4}e^{-d^2/4t}\sum_{i=0}^{k}t^i\Phi_i}_{C^{\ell}}\leq C_{j,k,\ell} t^j,
\]
where $||\cdot||_{C^{\ell}}$ denotes the $C^{\ell}$-norm of $C^l(P_{Spin}(M)\times P_{Spin}(M),\End(\wedge_{\bC}\bC^{n/2}))$. Moreover, the formal sum $\partial_t^k(4\pi t)^{-n(n+1)/4}e^{-d^2/4t}\sum_{j=0}^{\infty}t^j\Phi_j$ remains asymptotic to $\partial_t^kh_t$ under any finite amount of partial derivatives in $t$.
\end{theorem}
Applying \eqref{hka} to our equivariant formula for $k_t$, the compactness of $\Spin(n)$ shows that the equivariant averages remain asymptotic. Thus, we have
\[k_t(p,q)\sim\sum_{j=0}^{\infty}t^j(4\pi t)^{-n(n+1)/4}\int_{\Spin(n)}e^{-d(p\cdot g,q)^2/4t}\Phi_j(p\cdot g,q)\rho(g)dg
\]
\par
\subsection{The Getzler Rescaling on $P_{Spin}(M)$}
Let us now discuss what the Getzler rescaling looks like on $P_{Spin}(M)$. Consider $\kappa\in\Gamma(\mathscr{S}\boxtimes\mathscr{S}^*)$, $y\in M$, $p\in P_{Spin}(M)$ with $\pi(p)=y$. Note that under the isomorphism $\mathcal{H}_p\simeq T_yM$, if $v\in\mathcal{H}_p$, one has $\pi \exp^S_p(v)=\exp^M_y(v)$. Thus if we think of $\exp^S_p(v)$ as the radial parallel propagation of the spin frame $p$ along $\exp^M_y(v)$, then we have that the operator $P_{\exp^M_yv\rightarrow y}\kappa(\exp^M_y(v),y)$ on $\mathscr{S}_y$ has matrix $\kappa(\exp^S_p(v),p)\in\End(\wedge\bC^{n/2})$ if we use the spin frame $p$ at $y$. Combining this with $\End(\wedge\bC^{n/2})\simeq_c \bC\ell((\bR^n)^*)\simeq_{\sigma}\wedge (\bR^n)^*$, the matrix of the Getzler rescaled symbol $\delta_{u,y}(\kappa)(v)$ relative to the spin frame $p$ is
\[\sum_{j=0}^nu^{-j}\sigma_j(\kappa(\exp^S_p(uv),p)),\hspace{.5cm} 0<|u|\leq 1.
\]
Applying the parabolic Getzler rescaling to our heat kernel asymptotic expansion above, we see that the asymptotic sum of interest is
\begin{align*}
    \delta_{u,y}&(k_t)(v)=\sum_{j=0}^nu^{-j}\sigma_j(k_{u^2t}(\exp^S_p(uv),p))\\
    &\sim \sum_{k=0}^{\infty}\sum_{j=0}^n(u^2t)^k(4\pi u^2t)^{-n(n+1)/4}u^{-j}\int_{\Spin(n)}e^{-d(\exp^S_p(uv)\cdot g,p)^2/4u^2t}\Phi_k(\exp^S_p(uv)\cdot g,p)\sigma_j(\rho(g))dg\\
\end{align*}
We will focus for now on a particular term $\Phi_j$ in the asymptotic expansion above, which gives
\[(4\pi u^2t)^{-n(n+1)/2}\sum_{j=0}^nu^{-j}\int_{\Spin(n)}e^{-d(\exp^S_p(uv)\cdot g,p)^2/4u^2t}\Phi_j(\exp^S_p(uv)\cdot g,p)\sigma_j(\rho(g))dg
\]
The last simplification we make is to note that we can assume (by compactness of $P_{Spin}(M)$) that the neighborhood $U_{\Delta}$ is contained in an exponential tubular neighborhood of the diagonal $\Delta_{P_{Spin}}\subset P_{Spin}(M)\times P_{Spin}(M)$. Thus, we may replace the above integral over $\Spin(n)$ with an integral over $\mathfrak{spin}(n)$, and obtain
\[(4\pi u^2t)^{-n(n+1)/2}\sum_{j=0}^nu^{-j}\int_{\mathfrak{spin}(n)}e^{-d(\exp^S_p(uv)\cdot \exp A,p)^2/4u^2t}\Phi_j(\exp^S_p(uv)\cdot \exp A,p)\sigma_j(\rho(\exp A))J(A)dA
\]
where $J(A)$ is the Jacobian factor for exponential coordinates.
It will become apparent later that each such term in this asymptotic expansion has parabolic Getzler order $n$. Since each such $\Phi_j$ term has a factor of $(u^2t)^j$, we have that the only term that will be of any importance when computing the Getzler symbol of $k_t$ is $\Phi_0$.\par
\subsection{The Geometry of $P_{Spin}(M)$}
We let $\Omega\in\Omega^2_{Ad,h}(P_{Spin}(M),\mathfrak{spin}(n))$ denote the curvature of the distribution $\mathcal{H}$. Given $p\in P_{Spin}(M)$, $A\in\mathfrak{spin}(n)$, we define the operator $\tau(A\cdot \Omega_p)\in\mathfrak{so}(\mathcal{H}_p)$ via 
\[\langle \tau(A\cdot\Omega_p)v,w\rangle_{\mathcal{H}_p} :=2\langle A, \Omega_p(v,w)\rangle_{\mathfrak{spin}(n)}
\]
where $v,w\in\mathcal{H}_p$. Consider $x\in M$ and a point $p\in P_{Spin}(M)_x$. We have 
\[ \exp^S_p:T_pP_{Spin}(M)\rightarrow P_{Spin}(M)
\]
A bit of Jacobi field gymnastics (see \cite{BGV} Ch.5 for details) will show that given $A\in\mathfrak{spin}(n)$ small, and using $T_qP_{Spin}(M)\simeq\mathfrak{spin}(n)\oplus T_xM$ for $\pi(q)=x$, we have
\[(\exp^S_p)_{*,A}=\left[
\begin{array}{cc}
\frac{I-\exp(-ad_A)}{ad_A} & 0  \\ 
0 & \frac{I-\exp(-\tau(A\cdot\Omega_p)/2)}{\tau(A\cdot\Omega_p)/2}
\end{array}\right]
\]
In particular
\[\text{det}(\exp^S_p)_{*,A}=\text{det}\left(\frac{I-\exp(-ad_A)}{ad_A}\right)\cdot\text{det}\left(\frac{I-\exp(-\tau(A\cdot\Omega_p)/2)}{\tau(A\cdot\Omega_p)/2}\right)
\]
Thus, when computing $\Phi_0(p\cdot\exp A,p)$, we have
\begin{align*}
    \Phi_0(p\cdot \exp A,p)&=\text{det}^{-1/2}(\exp^S_p)_{*,A}\\
    &=\text{det}^{-1/2}\left(\frac{\sinh ad_A/2}{ad_A/2}\right)\cdot\text{det}^{-1/2}\left(\frac{\sinh\tau(A\cdot\Omega_p)/4)}{\tau(A\cdot\Omega_p)/4}\right).
\end{align*}
Next, consider $p\in P_{Spin}(M)$, $v\in\mathcal{H}_p$ and $A\in\mathfrak{spin}(n)$, both small. Then further Jacobi field backflips (see \cite{BGV} Ch.6) lead one to the small $u$ Taylor expansion
\[d(\exp^S_puv\cdot \exp A,p)^2= |A|^2+u^2\Big\langle v,\left(\frac{I-\exp(-\tau(A\cdot\Omega_p)/2)}{\tau(A\cdot\Omega_p)/2}\right)^{-1}v\Big\rangle + u^3f(u,A)
\] 
for a smooth function $f$ in $(u,A)$. Now observe that since $\tau(A\cdot\Omega_p)\in\mathfrak{so}(\mathcal{H}_p)$, we have $\langle v,\tau(A\cdot\Omega_p)v\rangle=0$. Thus we have 
\[\Big\langle v,\left(\frac{I-\exp(-\tau(A\cdot\Omega_p)/2)}{\tau(A\cdot\Omega_p)/2}\right)^{-1}v\Big\rangle=\langle v,(\tau(A\cdot\Omega_p)/4)\coth(\tau(A\cdot\Omega_p)/4)v \rangle
\]
\subsection{The Grassmann Calculus}
We recall (\cite{BGV} Ch.5) how one is to compute small time asymptotics of integrals of the form
\[(4\pi t)^{-n(n-1)/4}\int_{\mathfrak{spin}(n)}e^{-|A|^2/4t}\Psi(A)\exp_C(A)dA
\]
where $\exp_C$ denotes the exponential of $A$
in $\bC\ell((\bR^n)^*)$ and $\Psi$ is a decent smooth function on $\mathfrak{spin}(n)$. The key is to define what it means for $\Psi$ to be evaluated on a differential form in $\wedge^{even}(\bR^n)^*$. To do this, note that if we restrict $\sigma:\bC\ell((\bR^n)^*)\xrightarrow{\sim}\wedge(\bR^n)^*$ to $\mathfrak{spin}(n)$, we get the isomorphism 
\[\mathfrak{spin}(n)\simeq\wedge^2 (\bR^n)^*, \hspace{1cm} \sigma(e^ie^j)=e^i\wedge e^j
\]
Given $A=\sum_{i<j}a_{ij}e^ie^j\in\mathfrak{spin}(n)$ and $\Psi\in C^{\infty}(\mathfrak{spin}(n))$, a Taylor expansion at $A$ gives an element
\[
    T_A(\Psi)\in S^{\infty}(\mathfrak{spin}(n)^*),
\]
\[\Psi(A+x)\sim\sum_{\alpha}\frac{1}{\alpha !}\big(\partial_{e^1e^2}^{\alpha_{12}}\cdots\partial_{e^{n-1}e^n}^{\alpha_{n-1,n}}\Psi\big)(A)\prod_{i<j}x_{ij}^{\alpha_{ij}}
\]
Taking $e^ie^j$ as orthonormal gives an Ad-$\Spin(n)$-invariant inner product on $\mathfrak{spin}$ which allows us to identify $\mathfrak{spin}(n)$ with its dual $\mathfrak{spin}(n)^*$. The evaluation of $\Psi$ at $A+e\wedge e^*$ is defined by composing the Taylor expansion $T_A(\Psi)$ with $\sigma:S^{\infty}(\mathfrak{spin}(n)^*)\rightarrow\wedge^{even}(\bR^n)^*$
\[
    \Psi(A+e\wedge e^*):=\sigma\left( T_A(\Psi)\right)\in\wedge^{even}(\bR^n)^*,
\]
\[\Psi(A+e\wedge e^*)=\sum_{\alpha}\frac{1}{\alpha !}\big(\partial_{e^1e^2}^{\alpha_{12}}\cdots\partial_{e^{n-1}e^n}^{\alpha_{n-1,n}}\Psi\big)(A)\prod_{i<j}(e^i\wedge e^j)^{\alpha_{ij}}\in\wedge^{even}(\bR^n)^*
\]
To give two examples demonstrating the utility of the Grassmann calculus, we first cite Berline-Vergne who showed that, when evaluating $\Phi_0(p\cdot\exp A,p)$ on $2e\wedge e^*$ at $A=0$, one recovers 
\[ \Phi_0(2e\wedge e^*)=\text{det}^{1/2}(\frac{R_y/2}{\sinh R_y/2}),
\]
where $\pi(p)=y$.
Our extension is to use the Grassman calculus to compute the second order term of $d(\exp_p(uv)\cdot\exp A,p)^2$ on $2e\wedge e^*$. What we find is that 
\[\langle v,(\tau(A\cdot\Omega)/4)\coth(\tau(A\cdot\Omega/4)v\rangle(2e\wedge e^*)=\langle v,\frac{R_y}{2}\coth(\frac{R_y}{2})v\rangle.
\]
Turning back to the evaluation of small time asymptotics of Gaussian integrals on $\mathfrak{spin}(n)$, consider $\Psi$ a decent smooth function on $\mathfrak{spin}(n)$. Recall that when $A$ is close to zero, when evaluating $\sigma:\bC\ell((\mathbb{R}^n)^*)\xrightarrow{\sim}\wedge((\mathbb{R}^n)^*)^{\mathbb{C}}$ on $\exp_CA$, we have
\[
    \sigma(\exp_CA)=\mathscr{H}(\tau(A))\left(\exp_{\wedge}\sigma(A)\right),
\]
\[
    \mathscr{H}(\tau(A))=\text{det}^{1/2}\left(\frac{\sinh(\tau(A)/2)}{\tau(A)/2}\right)\text{det}^{-1/2}\left(\frac{\tanh(\tau(A)/2)}{\tau(A)/2}\right)\bigwedge\left(\frac{\tanh(\tau(A)/2)}{\tau(A)/2}\right)^{1/2}.
\]
Note that $\mathscr{H}(\tau(A))$ is an \emph{operator} on $\wedge(\bR^n)^*$. Thus, when apply the rescaling $t^{-j/2}$ on the degree $j$ component of $\sigma(\exp_CA)$, denoted $\delta_t\sigma(\exp_CA)$, we get
\[
    \delta_t\sigma(\exp_CA)=\mathscr{H}(\tau(A))(\exp_{\wedge}\sigma(A)/t).
\]
Applying this to our Gaussian integral gives
\begin{align*}
    (4\pi t)^{-n(n-1)/4}\int_{\mathfrak{spin}(n)}e^{-|A|^2/4t}\Psi(A)\delta_t\sigma(\exp_CA)dA&=\\
    (4\pi t)^{-n(n-1)/4}\int_{\mathfrak{spin}(n)}\mathscr{H}(\tau(A))&\exp\left(-\frac{|A-2e\wedge e^*|^2}{4t}\right)\Psi(A)dA
\end{align*}
which converges as $t\rightarrow 0^+$ to $\Psi(2e\wedge e^*)$.
It is clear that if $\Psi$ also depends smoothly on $t$, for $t\in\bR$, then a Taylor expansion in $t$ allows for the following small $t$ asymptotic expansion
\[
    (4\pi t)^{-n(n-1)/4}\int_{\mathfrak{spin}(n)}e^{-|A|^2/4t}\Psi(t,A)\delta_t\sigma(\exp_CA)dA\sim\sum_{k=0}^{\infty}\frac{t^k}{k!}\partial_t^k\Psi(0,2e\wedge e^*).
\]
In particular,
\[
    \lim_{t\rightarrow 0^+}(4\pi t)^{-n(n-1)/4}\int_{\mathfrak{spin}(n)}e^{-|A|^2/4t}\Psi(t,A)\delta_t\sigma(\exp_CA)dA=\Psi(0,2e\wedge e^*).
\]
To hint at the appearance of Mehler's kernel, we have that if
\[
    \Psi(v,A)=\exp(-\frac{1}{4}\langle v,(\tau(A\cdot\Omega)/4)\coth(\tau(A\cdot\Omega/4)v\rangle)\Phi_0(p\cdot\exp A,p)
\]
then evaluation of the Gaussian-Grassmann integral gives
\begin{align*}
(4\pi t)^{-n(n-1)/4}\int_{\mathfrak{spin}(n)}\mathscr{H}(\tau(A))&\exp\left(-\frac{|A-2e\wedge e^*|^2}{4t}\right)\Psi(v,A)dA\\
&\xrightarrow[t\rightarrow 0^+]{}\text{det}^{1/2}\left( \frac{R_y/2}{\text{sinh}(R_y/2)}\right)\text{exp}\left( -\frac{1}{4}\left\langle v\middle|\frac{R_y}{2}\text{coth}\left( \frac{R_y}{2}\right)\middle|v\right\rangle\right).
\end{align*}
The key to extracting the Mehler kernel from the Getzler rescaled heat kernel asymptotic expansion will be to reduce the analysis to the above simple computation.
\section{The Getzler Symbol of The Heat Kernel}
We begin our analysis of the Getzler rescaled heat kernel by restricting our attention to a single term of our asymptotic expansion and show how a Mehler-type kernel appears after rescaling. Thus consider the following term
\[(4\pi u^2t)^{-n(n+1)/4}\sum_{j=0}^nu^{-j}\int_{\mathfrak{spin}(n)}e^{-d(\exp^S_p(uv)\cdot \exp A,p)^2/4u^2t}\Phi(\exp^S_p(uv)\cdot \exp A,p)\sigma_j(\rho(\exp A))J(A)dA,
\]
where $\Phi$ is a smooth function on $P_{Spin}(M)\times P_{Spin}(M)$ with support contained in $V_{\Delta}$. Using the small time asymptotics of the distance function on $P_{Spin}(M)$, we can write the above integral as
\[(4\pi u^2t)^{-n(n+1)/4}\sum_{j=0}^{n/2}u^{-2j}\int_{\mathfrak{spin}(n)}e^{-|A|^2/4u^2t} \Psi(u,t,A)\sigma_{2j}(\exp_CA)dA
\]
with 
\[\Psi(u,t,v,A)=e^{-\langle v,(\tau(A\cdot\Omega_p)/4)\coth(\tau(A\cdot\Omega_p)/4)v \rangle/4t} e^{-uf(u,A)/4t}\Phi(\exp_p^S(uv)\cdot\exp A,p)J(A)
\]
Note that when $u=0$, we have 
\[\Psi(0,t,v,A)=e^{-\langle v,(\tau(A\cdot\Omega_p)/4)\coth(\tau(A\cdot\Omega_p)/4)v \rangle/4t}\Phi(p\cdot\exp A,p)J(A).
\]
Before doing anything too hasty, let us first re-express the above integral in a simpler form and rescale the integral by $\sqrt{t}$, thus giving:
\[
    (4\pi t)^{-n/2}u^{-n}\int_{\mathfrak{spin}(n)}\mathscr{H}(\sqrt{t}\tau(A))\exp\left(-\frac{|A-2\sqrt{t}e\wedge e^*|^2}{4u^2}\right)\Psi(u,t,v,\sqrt{t}A)dA.
\]
Recall that when computing the Getzler symbol we are concerned with the behaviour as $u\rightarrow 0$ with fixed $(t,v)$. Multiplying the above integral by $u^n$, we see that we have convergence as $u\rightarrow 0$ to the Mehler-type kernel
\[
    (4\pi t)^{-n/2}\Psi(0,t,v,2te\wedge e^*)=(4\pi t)^{-n/2}\Phi(2te\wedge e^*)\text{exp}\left( -\frac{1}{4}\left\langle v\middle|\frac{tR_y}{2}\text{coth}\left( \frac{tR_y}{2}\right)\middle|v\right\rangle\right).
\]
Thus we have that the parabolic Getzler rescaling, wehn multiplied by $u^n$, extends continuously to $u=0$. To see why the extension is \emph{smooth}, simply note that
\[
    \partial_u(4\pi u^2)^{-n(n-1)/4}\exp\left(-\frac{|A-2\sqrt{t}e\wedge e^*|^2}{4u^2}\right)=2u\Delta_{\mathfrak{spin}(n)}(4\pi u^2)^{-n(n-1)/4}\exp\left(-\frac{|A-2\sqrt{t}e\wedge e^*|^2}{4u^2}\right),
\]
and using integration by parts we may move the flat Laplacian $\Delta_{\mathfrak{spin}(n)}$ to the other terms within the Gaussian-Grassmann integral. Thus, we have shown that \emph{every} term in the equivariantly averaged asymptotic expansion has parabolic Getzler order $n$. However, the $j$-th term in the asymptotic expansion carries a $(u^2t)^j$ term which annihilates the smooth extension at $u=0$. Thus, the only term within the asymptotic expansion which extends smoothly to $u=0$, and is \emph{not} zero, is the $j=0$-term whose smooth extenstion to $u=0$ is
\[
    (4\pi t)^{-n/2}\text{det}^{1/2}\left( \frac{tR_y/2}{\text{sinh}(tR_y/2)}\right)\text{exp}\left( -\frac{1}{4t}\left\langle v\middle|\frac{tR_y}{2}\text{coth}\left( \frac{tR_y}{2}\right)\middle|v\right\rangle\right)
\]
Lastly, the $C^k$ estimates on the equivariantly averaged remainder term of any sufficiently long finite expansion show that the remainder term will extend to $0$ at $u=0$ in a differentiable manner up to order $k$, for $k$ as large as we desire. Thus, we have shown that the parabolic Getzler order of the spinor heat kernel $k_t$ is $n$ with
\[
    \lim_{u\rightarrow 0}u^{n}\delta_{u,y}(k_t)(v)=(4\pi t)^{-n/2}\text{det}^{1/2}\left( \frac{tR_y/2}{\text{sinh}(tR_y/2)}\right)\text{exp}\left( -\frac{1}{4t}\left\langle v\middle|\frac{tR_y}{2}\text{coth}\left( \frac{tR_y}{2}\right)\middle|v\right\rangle\right)
\]

\end{document}